# Ramanujan's q-continued fractions and Schröder-like numbers

## Johann Cigler

Fakultät für Mathematik
Universität Wien
A-1090 Wien, Nordbergstraße 15

johann.cigler@univie.ac.at

**Abstract**

In a recent paper G. Bhatnagar has given simple proofs of some of Ramanujan's continued fractions. In this note we show that some variants of these continued fractions are generating functions of $q-$ Schröder-like numbers.

**1. Introduction**

In a recent "tutorial" Gaurav Bhatnagar [3] has given simple proofs of some of Ramanujan's (convergent) $q-$ continued fractions by using an elementary method of Euler. I had not been aware of these continued fractions before but came across similar formulae in the study of formal power series which are generating functions of Schröder-like numbers and their $q-$ analogues (cf. [4]). The purpose of this note is to call attention to these connections and to give simple proofs of the corresponding continued fractions from this point of view.

A well-known example of the following story is the sequence of little Schröder numbers $(s_n)_{n \geq 0} = (1,1,3,11,45,197,\cdots)$ (cf. [6], A001003) whose generating function

$$f(z) = \sum_{n \geq 0} s_n z^n \tag{1.1}$$

satisfies

$$f(z) = 1 - zf(z) + 2zf(z)^2. \tag{1.2}$$

The computation of Hankel determinants for Schröder numbers leads to the continued fraction

$$f(z) = \cfrac{1}{1-z-}\cfrac{2z^2}{1-3z-}\cfrac{2z^2}{1-3z-}\cfrac{2z^2}{1-3z-}\cdots \tag{1.3}$$

but there are also other interesting continued fractions for the little Schröder numbers (cf. [6], A001003):

$$f(z) = \cfrac{1}{1+z-}\cfrac{2z}{1+z-}\cfrac{2z}{1+z-}\cfrac{2z}{1+z-}\cdots, \tag{1.4}$$

$$f(z) = \frac{1}{1-} \frac{z}{1-z-} \frac{z}{1-z-} \frac{z}{1-z-} \cdots \qquad (1.5)$$

and

$$f(z) = \frac{1}{1-} \frac{z}{1-} \frac{2z}{1-} \frac{z}{1-} \frac{2z}{1-} \frac{z}{1-} \cdots. \qquad (1.6)$$

These will appear as special cases of the following considerations.

## 2. Generating functions of q-Schröder-like numbers

Some of the following results have been obtained in [4]. We repeat them in order to make the exposition self-contained.
Let $x, y$ be real or complex numbers and $z$ an indeterminate.

Let the formal power series

$$F(z) = F(z, x, y) = \sum_{n \geq 0} A(n, x, y) z^n \qquad (2.1)$$

satisfy the identity

$$F(z) = 1 + xzF(z) + yzF(z)F(qz). \qquad (2.2)$$

This implies that

$$A(n, x, y) = xA(n-1, x, y) + y \sum_{k=0}^{n-1} A(k, x, y) q^k A(n-1-k, x, y) \text{ with } A(0, x, y) = 1.$$

We are mainly interested in the series

$$f(z) = f(z, x, y) = \sum_{n \geq 0} a(n, x, y) z^n = \frac{x + yF(z, x, y)}{x + y}. \qquad (2.3)$$

It is easily verified that it satisfies the equation

$$f(z, x, y) = 1 - xzf(qz, x, y) + (x+y)zf(z, x, y)f(qz, x, y). \qquad (2.4)$$

Its coefficients are given by

$$a(n, x, y) = -q^{n-1}xa(n-1, x, y) + (x+y)\sum_{k=0}^{n-1} a(k, x, y) q^k a(n-1-k, x, y) \text{ with } a(0, x, y) = 1.$$

The sequence $(A(n,1,q))$ is a $q-$ analogue of the (large) Schröder numbers and the sequence $(a(n,1,q))$ is a $q-$ analogue of the little Schröder numbers. The numbers $A(n,1,q)$ have been studied from a combinatorial point of view in [1]. We call $A(n, x, y)$ and $a(n, x, y)$ $q-$ Schröder-like numbers. They are polynomials in $x$ and $y$.
The numbers $A(n,0,1) = a(n,0,1) = C_n(q)$ are the Carlitz $q-$ Catalan numbers. For $q = 1$ they reduce to the Catalan numbers $C_n = \frac{1}{n+1}\binom{2n}{n}$.

Equation (2.2) implies immediately the expansion of the formal power series $F(z,x,y)$ into a continued fraction

$$F(z,x,y) = \frac{1}{1-xz-yzF(qz,x,y)} = \cfrac{1}{1-xz-\cfrac{yz}{1-qxz-\cfrac{qyz}{1-q^2xz-\cdots}}}. \qquad (2.5)$$

Using (2.3) it is easily verified that

$$F(z,x,y) = 1 + (x+y)zF(z,x,y)f(qz,x,y). \qquad (2.6)$$

Consider the uniquely determined series $h(z) = h(z,x,y) = 1 + \sum_{n\geq 1} h_n(x,y)z^n$ which satisfies

$$F(z,x,y) = \frac{h(qz,x,y)}{h(z,x,y)}. \qquad (2.7)$$

From the defining equation for $F(z)$ we get $\dfrac{h(qz)}{h(z)} = 1 + xz\dfrac{h(qz)}{h(z)} + yz\dfrac{h(qz)}{h(z)}\dfrac{h(q^2z)}{h(qz)}$

and therefore $h(qz) = h(z) + xzh(qz) + yzh(q^2z)$.

Comparing coefficients we get
$(q^n - 1)h_n = q^{n-1}\left(x + q^{n-1}y\right)h_{n-1}$.
This implies

$$h_k(x,y) = q^{\binom{k}{2}} \frac{(x+y)(x+qy)\cdots(x+q^{k-1}y)}{(q-1)(q^2-1)\cdots(q^k-1)} \qquad (2.8)$$

and thus

$$h(z,x,y) = \sum_{k\geq 0} q^{\binom{k}{2}} \frac{(x+y)(x+qy)\cdots(x+q^{k-1}y)}{(1-q)(1-q^2)\cdots(1-q^k)}(-z)^k. \qquad (2.9)$$

The identity

$(x+y)h(z,x,qy) = xh(z,x,y) + yh(qz,x,y),$

implies

$$f(z,x,y) = \frac{x + yF(z,x,y)}{x+y} = \frac{h(z,x,qy)}{h(z,x,y)}. \qquad (2.10)$$

Observing that

$$F(z,x,y)f(qz,x,y) = \frac{h(qz,x,y)}{h(z,x,y)} \frac{h(qz,x,qy)}{h(qz,x,y)} = \frac{h(qz,x,qy)}{h(z,x,y)}$$
$$= \frac{h(z,x,qy)}{h(z,x,y)} \frac{h(qz,x,qy)}{h(z,x,qy)} = f(z,x,y)F(z,x,qy)$$
(2.11)

we see that (2.4) and (2.6) can be written in the form

$$F(z,x,y) = 1 + (x+y)zF(z,x,y)f(qz,x,y) \tag{2.12}$$

and

$$f(z,x,y) = 1 + yzf(z,x,y)F(z,x,qy). \tag{2.13}$$

The last line follows from (2.3) and (2.11):

$$f(z,x,y) = \frac{x + yF(z,x,y)}{x+y} = \frac{x+y+y(F(z,x,y)-1)}{x+y} = 1 + yzF(z,x,y)f(qz,x,y).$$

We shall also need the formula

$$h(z,x,y) = \left(1 + \frac{xz}{q}\right)h(z,x,qy) - \frac{z}{q}(x+qy)h\left(z,x,q^2 y\right). \tag{2.14}$$

For the proof it suffices to show that

$$h_k(x,y) = (-1)^k q^{\binom{k}{2}} \frac{(x+y)(x+qy)\cdots(x+q^{k-1}y)}{(1-q)(1-q^2)\cdots(1-q^k)} \tag{2.15}$$

satisfies

$$h_k(x,y) = h_k(x,qy) + \frac{x}{q}h_{k-1}(x,qy) - (x+qy)\frac{1}{q}h_{k-1}(x,q^2 y). \tag{2.16}$$

This is equivalent with

$$q^{\binom{k}{2}} \frac{(x+y)(x+qy)\cdots(x+q^{k-1}y)}{(1-q)(1-q^2)\cdots(1-q^k)} = q^{\binom{k}{2}} \frac{(x+qy)(x+qy)\cdots(x+q^k y)}{(1-q)(1-q^2)\cdots(1-q^k)}$$
$$-q^{\binom{k-1}{2}} \frac{x}{q} \frac{(x+qy)(x+q^2 y)\cdots(x+q^{k-1}y)}{(1-q)(1-q^2)\cdots(1-q^{k-1})} + (x+qy)q^{\binom{k-1}{2}-1} \frac{(x+q^2 y)(x+q^3 y)\cdots(x+q^k y)}{(1-q)(1-q^2)\cdots(1-q^{k-1})}$$

or equivalently
$$q^k(x+y) = q^k(x+q^k y) + x(q^k-1) - (x+q^k y)(q^k-1),$$

which is obviously true.

Let now $y$ also be an indeterminate. Then we can give another characterization of $f(z,x,y)$:

$$f(z,x,y) = \frac{\sum_{k\geq 0} \frac{q^{k^2}}{(q;q)_k (xz;q)_k}(-yz)^k}{\sum_{k\geq 0} \frac{q^{k^2-k}}{(q;q)_k (xz;q)_k}(-yz)^k}. \tag{2.17}$$

Here as usual $(xz;q)_k = (1-xz)(1-qxz)\cdots(1-q^{k-1}xz)$.

To prove this observe that (2.14) implies

$$f(z,x,y) = \frac{h(z,x,qy)}{h(z,x,y)} = \frac{1}{1+\frac{xz}{q}-\frac{z}{q}(x+qy)\frac{h(z,x,q^2y)}{h(z,x,qy)}} = \frac{1}{1+\frac{xz}{q}-\frac{z}{q}(x+qy)f(z,x,qy)} \tag{2.18}$$

and thus also

$$f(z,x,y) = 1 - \frac{xz}{q}f(z,x,y) + \frac{z}{q}(x+qy)f(z,x,y)f(z,x,qy). \tag{2.19}$$

If we set in this equation

$$f(z,x,y) = \frac{H(z,x,qy)}{H(z,x,y)} \tag{2.20}$$

with a formal power series

$$H(z,x,y) = \sum_n H_n(x,z) y^n \tag{2.21}$$

(2.19) implies

$$H(z,x,qy) = H(z,x,y) - \frac{xz}{q}H(z,x,qy) + \frac{z}{q}(x+qy)H(z,x,q^2y). \tag{2.22}$$

Comparing coefficients of $y^n$ we get

$$q^n H_n = H_n - q^{n-1}xzH_n + q^{2n-1}xzH_n + q^{2n-2}zH_{n-1}$$

or

$$H_n(x,z) = -\frac{q^{2n-2}z}{(1-q^n)(1-q^{n-1}xz)}H_{n-1}(x,z).$$

This gives

$$H(z,x,y) = \sum_{k\geq 0} \frac{q^{k^2-k}}{(q;q)_k (xz;q)_k}(-yz)^k \qquad (2.23)$$

as a formal power series in $y$ whose coefficients are formal power series in $z$.

Since $\dfrac{1}{(xz;q)_k} = \sum_{j\geq 0}\begin{bmatrix} k+j-1 \\ j \end{bmatrix} z^j$ we get $H(z,x,y) = \sum_{n\geq 0} z^n \sum_{j+k=n}\begin{bmatrix} n-1 \\ j \end{bmatrix}(-1)^k \dfrac{q^{k^2-k}}{(q;q)_k} y^k.$

Here $\begin{bmatrix} n \\ k \end{bmatrix} = \begin{bmatrix} n \\ k \end{bmatrix}_q = \dfrac{(q;q)_n}{(q;q)_k (q;q)_{n-k}}$ is a $q$-binomial coefficient.

Thus $H(z,x,y)$ is a formal power series in $z$ whose coefficients are polynomials in $y$. Therefore the right-hand side of (2.17) is also a formal power series in $z$ whose coefficients are polynomials in $y$. Therefore it is possible in (2.17) to replace the indeterminate $y$ by a real or complex number.

Comparing (2.20) with (2.10) we see that

$$f(z,x,y) = \frac{h(z,x,qy)}{h(z,x,y)} = \frac{H(z,x,qy)}{H(z,x,y)}. \qquad (2.24)$$

This implies

$$\frac{h(z,x,qy)}{H(z,x,qy)} = \frac{h(z,x,y)}{H(z,x,y)}. \qquad (2.25)$$

Since $\dfrac{h(z,x,y)}{H(z,x,y)} = \sum_n c_n(x,z) y^n$ is a formal power series in $y$ whose coefficients are formal power series in $z$ the equation $\sum_n c_n(x,z) y^n = \sum_n c_n(x,z) q^n y^n$ implies $c_n(x,z) = 0$ for $n > 0$.

Thus

$$\frac{h(z,x,y)}{H(z,x,y)} = \frac{h(z,x,0)}{H(z,x,0)} = \frac{1}{e(xz)}. \qquad (2.26)$$

Here $e(z) = \sum_{n\geq 0} \dfrac{z^n}{(q;q)_n}$ denotes the $q$-exponential series which satisfies $\dfrac{1}{e(z)} = \sum_{n\geq 0} (-1)^n q^{\binom{n}{2}} \dfrac{z^n}{(q;q)_n}.$

Formula (2.26) is a formal power series version of [2], Entry 9.

For $(x,y) = (q,-q)$ we have $h(z,q,-q) = 1$ and $H(z,q,-q) = \sum_{k\geq 0} \dfrac{q^{k^2} z^k}{(q;q)_k (qz;q)_k}.$

In this case (2.26) reduces to a well-known identity of Cauchy.

## 3. Associated continued fractions

From the results of [4], (3.20) (there is a typo; it should read $s(n) = q^{n-1}(x + q^n(1+q)y))$ we can deduce the Jacobi type continued fraction for $f(z,x,y)$ which we state without proof:

$$f(z,x,y) = \frac{1}{1-yz} - \frac{y(x+qy)z^2}{1-(x+q(1+q)y)z} - \frac{q^3 y(x+q^2 y)z^2}{1-q(x+q^2(1+q)y)z} - \frac{q^6 y(x+q^3 y)z^2}{1-q^2(x+q^3(1+q)y)z} - \cdots.$$

But here we are interested in other continued fractions.

**a)** From (2.12) and (2.13) we get

$$F(z,x,y) = \frac{1}{1-(x+y)zf(qz,x,y)} \tag{3.1}$$

and

$$f(z,x,y) = \frac{1}{1-yzF(z,x,qy)}. \tag{3.2}$$

This gives the following continued fraction for $f(z,x,y)$:

$$f(z,x,y) = \cfrac{1}{1 - \cfrac{yz}{1 - \cfrac{(x+qy)z}{1 - \cfrac{q^2 yz}{1 - \cfrac{q(x+q^2 y)z}{1 - \cfrac{q^4 yz}{1 - \cdots}}}}}} \tag{3.3}$$

which generalizes (1.6).

**Remark**

This and the other results about continued fractions are of course well known and due to Ramanujan who essentially developed the right-hand side of (2.24) into a convergent continued fraction of the form (3.3). The only new fact if anything in our approach is the connection with $q$-analogues of Schröder numbers. We are not interested in convergence questions and use only formal power series in $z$ instead of convergent power series in $q$. In this sense (3.3) can also be found in [2], (13.5) and [3], (6.5), where in our notation $f(z,qx,-qy)$ instead of $f(z,x,y)$ has been used.

**b)** Another continued fraction for $f(z)$ which is related to [3], (7.1) is

$$f(z,x,y) = \cfrac{1}{1+\cfrac{xz}{q} - \cfrac{\cfrac{z}{q}(x+qy)}{1+\cfrac{xz}{q} - \cfrac{\cfrac{z}{q}(x+q^2y)}{1+\cfrac{xz}{q} - \cdots}}} \tag{3.4}$$

which generalizes (1.4).

This is an immediate consequence of (2.18).

**Remark**

Note that (3.4) is essentially Touchard's continued fraction which has been studied by Helmut Prodinger in [5]. We get

$$f(z,q,-q) = \cfrac{1}{1+z - \cfrac{(1-q)z}{1+z - \cfrac{(1-q^2)z}{1+z - \cdots}}} \tag{3.5}$$

By (2.4) or by (2.10) we derive that

$$f(z,q,-q) = h(z,q,-q^2) = \sum_n (-1)^n q^{\binom{n+1}{2}} z^n. \tag{3.6}$$

Prodinger has given a direct proof that

$$\frac{H(z,q,-q^{i+2})}{H(z,q,-q)} = \sum_{k \geq 0} q^{\binom{k+1}{2}} \begin{bmatrix} k+i \\ k \end{bmatrix} (-z)^k. \tag{3.7}$$

In our setting this follows from (2.26) since

$$\frac{H(z,q,-q^{i+2})}{H(z,q,-q)} = \frac{h(z,q,-q^{i+2})}{h(z,q,-q)} = \sum_{k \geq 0} q^{\binom{k+1}{2}} \frac{(1-q^{i+1})(1-q^{i+2})\cdots(1-q^{i+k}y)}{(1-q)(1-q^2)\cdots(1-q^k)} (-z)^k.$$

**c)** Finally we derive the analogue of (1.5) (see [2], Entry 15, or [3], (5.3))

$$f(z,x,y) = \frac{1}{1-}\frac{yz}{1-xz-}\frac{qyz}{1-qxz-}\cdots. \tag{3.8}$$

By (2.12) we have

$$f(z,x,y) = \frac{1}{1-yzF(z,x,qy)}.$$

By (2.2)

$$F(z,x,qy) - xzF(z,x,qy) - qyzF(qz,x,qy)F(z,x,qy)$$

and therefore

$$F(z,x,qy) = \frac{1}{1-xz-qyzF(qz,x,qy)} = \frac{1}{1-xz-}\frac{qyz}{1-qxz-}\frac{q^2yz}{1-q^2xz-}\cdots.$$

This implies (3.8).

As a special case we get

$$\frac{1}{1+}\frac{z}{1-z+}\frac{qz}{1-qz+}\frac{q^2z}{1-q^2z+}\cdots = \sum_{n\geq 0}(-1)^n q^{\binom{n}{2}} z^n. \tag{3.9}$$

Note that for the famous Rogers-Ramanujan continued fraction

$$f(z,0,-q) = \frac{1}{1+}\frac{qz}{1+}\frac{q^2z}{1+}\frac{q^3z}{1+}\cdots \quad \text{both formulae (2.10) and (2.17) coincide.}$$

For the little $q$ – Schröder numbers the corresponding continued fractions are

$$f(z,1,q) = \frac{1}{1+z/q-}\frac{(1+q^2)z/q}{1+z/q}-\frac{(1+q^3)z/q}{1+z/q}-\cdots, \tag{3.10}$$

$$f(z,1,q) = \frac{1}{1-}\frac{qz}{1-z-}\frac{q^2z}{1-qz-}\cdots, \tag{3.11}$$

and

$$f(z,1,q) = \frac{1}{1-}\frac{qz}{1-}\frac{(1+q^2)z}{1}-\frac{q^3z}{1-}\frac{(q+q^4)z}{1}-\cdots. \tag{3.12}$$

**References**

[1] E. Barcucci, A. del Lungo, E. Pergola, and R. Pinzani, Some combinatorial interpretations of q-analogs of Schröder numbers, Ann. Comb. 3 (1999), 171 - 190

[2] B. C. Berndt, Ramanujan's Notebooks, Part III, Springer 1991

[3] G. Bhatnagar, How to prove Ramanujan's q-continued fractions, arXiv:1205.5455

[4] J. Cigler, Hankel determinants of Schröder-like numbers, arXiv: 0901.4680

[5] H. Prodinger, On Touchard's continued fraction and extensions: combinatorics-free, self-contained proofs,   http://math.sun.ac.za/~hproding/pdffiles/touchard-2011.pdf

[6] The On-Line Encyclopedia of Integer Sequences, http://oeis.org/